\documentclass[12pt]{article}

\usepackage{amsmath}
\usepackage{amsfonts}
\usepackage{amssymb}
\usepackage{bbm}

\newtheorem{theorem}{Theorem}
\newtheorem{lemma}{Lemma}
\newtheorem{corollary}{Corollary}
\newtheorem{assumption}{Assumption}

\newcommand{\R}{\mathbbm{R}}
\newcommand{\N}{\mathbbm{N}}
\newcommand{\C}{\mathbbm{C}}
\newcommand{\bfA}{\mathbf{A}}
\newcommand{\bfa}{\mathbf{a}}
\newcommand{\bfb}{\mathbf{b}}
\newcommand{\bfe}{\mathbf{e}}

\newcommand{\bfu}{\mathbf{u}}
\newcommand{\bfv}{\mathbf{v}}
\newcommand{\bfw}{\mathbf{w}}
\newcommand{\bfx}{\mathbf{x}}

\newcommand{\bfz}{\mathbf{z}}

\begin{document}

\title{A Randomized Parallel Algorithm with Run Time $O(n^2 )$ for Solving an $n \times n$ System of Linear Equations}

\author{J\"org Fliege\footnote{The author is indepted to Ian Hawke, School of Mathematics, University of Southampton,
                               for pointing out an error in a previous version of this note.} \\ %
%
School of Mathematics \\ %
The University of Southampton  \\ %
Southampton SO17 1BJ \\ %
UK \\ %
\texttt{J.Fliege@soton.ac.uk}  \\ %
Tel +44-23-8059-8453   \\ %
Fax +44-23-8059-3131    \\ %
\texttt{http://www.cormsis.soton.ac.uk}
}

\maketitle

In this note, following suggestions by Tao~\cite{Tao:2012},
we extend the randomized algorithm for linear equations over
prime fields by Raghavendra~\cite{Raghavendra:2012} to a randomized
algorithm for linear equations over the reals. We also show
that the algorithm can be parallelized to solve a system of linear 
equations $\bfA \bfx = \bfb$ with a regular $n \times n$ matrix~$\bfA$ 
in time $O(n^2)$, with probability one. Note that we do not
assume that~$\bfA$ is symmetric.

Let $m, n \in \N$ with $m \leq n$ and consider an $m \times n$ matrix  
$\bfA \in \R^{m \times n}$ as well as a right-hand side vector $\bfb \in \R^m$. 
There are many applications in which it is known in advance that~$\bfA$ has full
row rank, i.e.~the system of linear equations $\bfA \bfx = \bfb$ has at least one solution.
We are interested in solving such a system in the sense that we want to construct
a vector $\bfx \in \R^n$ that fulfils these equations, given the knowledge that
$\bfA$~has full row rank.
This has, of course, important applications for the case $m = n$.
For this problem, we consider the algorithm described below. 
In what follows, let $\bfA$ have the row vectors $\bfa_1, \ldots, \bfa_m \in \R^n$, i.e.\ 
\[
\bfA = 
\left[
\begin{array}{c}
\bfa_1^{\top} \\
\vdots \\
\bfa_m^{\top}
\end{array}
\right] .
\]

We consider a random vector $\xi : \Theta \longrightarrow \R^n$ defined on some
probability space $(\Theta, \mathcal{F}, P)$, for which the following holds.
\begin{assumption} \label{ass:1}
For arbitrary $\bfa \in \R^n$, $\bfa \neq 0$, and $\beta \in \R$ we have
\[
\mbox{\rm{Prob}} \left( \bfa^{\top} \xi = \beta \right)  = 0 .
\]
\end{assumption}
In other words, the random vector is not biased towards
particular affine subspaces of~$\R^n$. Examples for corresponding distributions include
the case in which each coordinate $\xi_i$ ($i = 1, \ldots, n$) is independently drawn
from a Gaussian distribution on~$\R$, or from a uniform distribution over a 
certain interval, 
or in which~$\xi$~is continuous uniformly distributed on the unit sphere
$\{ \bfx \in \R^n \mid  \sum_{i=1}^n x_i^2 = 1 \}$. From the assumption, it follows
readily that
\[
\mbox{\rm{Prob}} \left( \xi = \bfx \right)  = 0 .
\]
for all $\bfx \in \R^n$, as 
$\mbox{\rm{Prob}} \left( \xi = \bfx \right) 
 \leq \sum_{i = 1}^n \mbox{\rm{Prob}} \left( ( \bfe^{(i)} )^{\top} \xi = \bfx_i \right) = 0$,
where $\bfe^{(i)}$ are the Cartesian unit vectors ($i = 1, \ldots, n$). 
For technical reasons, we will also assume that $\xi(F)$ is a measurable set 
(in the usual sense of the natural Borel $\sigma$-Algebra of~$\R^n$) for all 
$F \in \mathcal{F}$.

We are now ready to state the main algorithm.

\begin{enumerate}
\item \textbf{Input:} $(\bfA, \bfb)$ with matrix $A$ as row vectors $\bfa_1, \ldots, \bfa_m \in \R^n$,
      and a right-hand side $\bfb \in \R^m$.
\item Let $\bfv_1, \bfv_2, \ldots, \bfv_{n+1} \in \R^n$ denote identically independent distributed 
      samples of the random variable~$\xi$.
\item \texttt{for} $k = 1, \ldots, m$ \texttt{do}
      \begin{enumerate}
      \item Choose $n+1$ random pairs $(i_1, j_1), \ldots, (i_{n+1}, j_{n+1})$ with
            $i_{\ell} < j_{\ell}$ for $\ell = 1, \ldots, n+1$, and all pairs unequal
            to each other.
      \item \texttt{for} $\ell = 1, \ldots, n+1$ \texttt{do} 
            \begin{enumerate}
            \item $\bfx_{\ell} := \mbox{\textbf{rec}}( \bfv_{i_{\ell}}, \bfv_{j_{\ell}}, \bfa_k, b_k )$
            \end{enumerate}
      \item \texttt{if} one of the calls to \mbox{\textbf{rec}} stops with failure,
            \texttt{then STOP with failure}
      \item Otherwise, set $\bfv_{\ell} := \bfx_{\ell}$ for $\ell = 1, \ldots, n$
      \end{enumerate}
\item \textbf{Output:} $\bfv_1, \ldots, \bfv_{n+1}$.
\end{enumerate}

This algorithm makes use of the subroutine \mbox{\textbf{rec}} ("recombination"), 
defined as follows:

\begin{enumerate}
\item \textbf{Input:} $(\bfu, \bfv, \bfa, \beta)$ with vectors 
      $\bfu, \bfv, \bfa \in \R^n$ and a real number~$\beta$.
\item \texttt{if} $\bfa^{\top} (\bfu-\bfv) = 0$ \texttt{then STOP with failure}
\item Otherwise, set
      \[
      t := \frac{\beta - \bfa^{\top} \bfv}{\bfa^{\top} (\bfu-\bfv)} 
      \]
      and set $\bfz := t \bfu + (1-t)\bfv$.
\item \textbf{Output:} $\bfz$.
\end{enumerate}

In what follows, we will show the following.

\begin{theorem} \label{theorem:main}
Suppose $\bfA$ has full row rank and that Assumption~\ref{ass:1} holds.
\begin{enumerate}
\item
With probability one, the randomized algorithm described above stops
after~$m$ steps with output $\bfv_{\ell} \in \R^n$,
$\ell = 1, \ldots, n+1$, such that  $\bfA \bfv_{\ell} = \bfb$ holds
for $\ell = 1, \ldots, n+1$.
\item
With probability one, the run time of the algorithm is bounded by
$O(n^2 m)$ floating point operations.
\end{enumerate}
\end{theorem}

From this, the following corollary immediately follows.
\begin{corollary}
Consider a regular matrix $\bfA \in \R^{n \times n}$, a right-hand side $\bfb \in \R^n$
and suppose that Assumption~\ref{ass:1} holds.
Then, with probability one, the randomized algorithm above solves the linear system 
of equations $\bfA \bfx = \bfb$ in $O(n^3 )$ floating point operations.
\end{corollary}

As it can be clearly seen, Step~3 of the algorithm can be fully parallelized.
As each call to \mbox{\textbf{rec}} costs $O(n)$ flops, we arrive at the
main result of this note.
\begin{corollary}
Consider a regular matrix $\bfA \in \R^{n \times n}$, a right-hand side $\bfb \in \R^n$,
and suppose that Assumption~\ref{ass:1} holds.
Then the randomized algorithm above can then be parallelized such that, with probability one, 
it solves the linear system  of equations $\bfA \bfx = \bfb$ in time $O(n^2 )$.
\end{corollary}

We start the analysis with a straightforward result.
\begin{lemma} \label{lemma:rec}
Consider vectors $\bfu, \bfv, \bfa \in \R^n$ with $\bfa \neq 0$ and a real 
number~$\beta$.
Then, either $\bfa^{\top} (\bfu-\bfv) = 0$ or the subroutine \mbox{\rm\textbf{rec}}
returns a vector $\bfz = t\bfu + (1-t)\bfv$ with $\bfa^{\top} \bfz = \beta$.
\end{lemma}
\textbf{Proof:} By construction.  $\Box$

Next, we consider the first~$i$ iterations of the algorithm.
\begin{lemma} \label{lemma:ind}
Let Assumption~\ref{ass:1} hold, 
let $1 \leq i \leq m$ and let $\bfa_1, \ldots, \bfa_i$ be linearly independent.
Then, the following holds.
\begin{enumerate}
\item With probability one the algorithm has not stopped with failure in the first~$i$
iterations and the vectors $\bfv_{\ell} = \bfv_{\ell}^{(i)}$ ($\ell = 1, \ldots, n+1$),
produced in step~$i$ of the algorithm, satisfy
$\bfa_j^{\top} \bfv_{\ell} = b_j$ 
for $j = 1, \ldots, i$ and $\ell = 1, \ldots, n+1$.
\item Let $\bfa \in \R^n$ be an arbitrary vector with $\bfa \neq 0$ 
      and let $\beta \in \R$ be arbitrary.
      Then, for all $\ell = 1, \ldots, n+1$,
      $\bfa^{\top} \bfv_{\ell} = \beta$ holds with probability zero,
      where $\bfv_{\ell} = \bfv_{\ell}^{(i)}$ denote the iteration vectors
      of the algorithm after step~$i$.
\end{enumerate}
\end{lemma}
\textbf{Proof.} We show both claims by induction.
\begin{enumerate}
\item $i =1$: claim~2 follows directly from Assumption~\ref{ass:1}. Claim~1
      follows from Lemma~\ref{lemma:rec}, as $\bfA$ has full rank and
      $\bfa_1^{\top} (\bfv_{i} - \bfv_j ) = 0$ holds with probability zero
      for all $i, j = 1, \ldots, n+1$, $i \neq j$.
\item $i \rightarrow i+1 \leq m$:
      We start with claim~1.
      Suppose that $\bfv_{\ell}$ satisfy $\bfa_j^{\top} \bfv_{\ell} = b_j$ 
      for $j = 1, \ldots, i$ and $\ell = 1, \ldots, n+1$.
      Let $(\ell_1, \ell_2)$ be a randomly chosen pair of indices with $1 \leq \ell_1 < \ell_2 \leq n+1$.
      Without loss of generality, assume $\ell_1 = 1$ and $\ell_2 = 2$.
      If $\mbox{\textbf{rec}} ( \bfv_1, \bfv_2, \bfa_{i+1}, b_{i+1} )$ returns
      a vector~$\bfx$ without failure, then $\bfa_{i+1}^{\top} \bfx = b_{i+1}$
      by Lemma~\ref{lemma:rec}. Also, $\bfx$ is a convex combination of~$\bfv_1$
      and~$\bfv_2$ and thus fulfils $\bfa_j^{\top} \bfx = b_j$ 
      for $j = 1, \ldots, i$. But, due to claim~2, the call 
      $\mbox{\textbf{rec}} ( \bfv_1, \bfv_2, \bfa_{i+1}, b_{i+1} )$
      returns without failure with probability one. This shows claim~1.
      
      It remains to perform the inductive step for claim~2. As above,
      let us choose the pair of vectors $\bfv_1 , \bfv_2$ without loss of
      generality.
      Due to the induction hypothesis, we have, with 
      probability one,
      \begin{eqnarray*}
      \bfx & = & \frac{b_i - \bfa_i^{\top} \bfv_2}{\bfa_i^{\top} (\bfv_1-\bfv_2)} \bfv_1
                   + \left( 1-
                     \frac{b_i - \bfa_i^{\top} \bfv_2}{\bfa_i^{\top} (\bfv_1-\bfv_2)} \right) \bfv_2 \\
      & = & \frac{1}{\bfa_i^{\top} (\bfv_1-\bfv_2)} 
            \left( 
            ( b_i - \bfa_i^{\top} \bfv_2 ) \bfv_1
            - ( b_i - \bfa_i^{\top} \bfv_1 ) \bfv_2
            \right)
      \end{eqnarray*}
      and therefore $\bfa^{\top} \bfx = \beta$ if and only if
      \[
      ( b_i - \bfa_i^{\top} \bfv_2 ) \bfa^{\top} \bfv_1 = 
      ( b_i - \bfa_i^{\top} \bfv_1 ) \bfa^{\top} \bfv_2 + \beta \bfa_i^{\top} (\bfv_1-\bfv_2) 
      \]
      holds. Thus,
      \begin{eqnarray*}
      \lefteqn{ \mbox{\rm{Prob}} \left( \bfa^{\top} \bfx = \beta \right) } \\
      & = & 
      \mbox{\rm{Prob}} \left( ( b_i - \bfa_i^{\top} \bfv_2 ) \bfa^{\top} \bfv_1 = 
      ( b_i - \bfa_i^{\top} \bfv_1 ) \bfa^{\top} \bfv_2 + \beta \bfa_i^{\top} (\bfv_1-\bfv_2) \right) \\
      & = & 
      \int_{-\infty}^{\infty}
      \mbox{\rm{Prob}} \left( ( b_i - \bfa_i^{\top} \bfv_2 ) \bfa^{\top} \bfv_1 = \zeta \mbox{~and~}
      \zeta = ( b_i - \bfa_i^{\top} \bfv_1 ) \bfa^{\top} \bfv_2 + \beta \bfa_i^{\top} (\bfv_1-\bfv_2) \right) 
      \, \mbox{d}\zeta \\
      & \leq & 
      \int_{-\infty}^{\infty}
      \mbox{\rm{Prob}} \left( ( b_i - \bfa_i^{\top} \bfv_2 ) \bfa^{\top} \bfv_1 = \zeta \right) 
      + 
      \mbox{\rm{Prob}} \left(  \zeta = ( b_i - \bfa_i^{\top} \bfv_1 ) \bfa^{\top} \bfv_2 + \beta \bfa_i^{\top} (\bfv_1-\bfv_2) \right) 
      \, \mbox{d}\zeta ,
      \end{eqnarray*}
      where the existence of the integrals are guaranteed as $\xi$ maps 
      measurable sets on measurable sets, by assumption.
      But
      \begin{eqnarray*}
      \lefteqn{
      \mbox{\rm{Prob}} \left(  \zeta = ( b_i - \bfa_i^{\top} \bfv_1 ) \bfa^{\top} \bfv_2 + \beta \bfa_i^{\top} (\bfv_1-\bfv_2) \right) 
      } \\
      & = & 
      \int_{-\infty}^{\infty}
      \mbox{\rm{Prob}} \left(  \zeta = ( b_i - \bfa_i^{\top} \bfv_1 ) \bfa^{\top} \bfv_2 + \beta \eta 
      \mbox{~and~} \eta =  \bfa_i^{\top} (\bfv_1-\bfv_2)
      \right)
      \, \mbox{d}\eta \\
      & \leq &
      \int_{-\infty}^{\infty}
      \mbox{\rm{Prob}} \left(  \zeta = ( b_i - \bfa_i^{\top} \bfv_1 ) \bfa^{\top} \bfv_2 + \beta \eta \right)
      + \mbox{\rm{Prob}} \left( \eta = \bfa_i^{\top} (\bfv_1-\bfv_2)
      \right) 
      \, \mbox{d}\eta \\
      & = & 
      \int_{-\infty}^{\infty}
      \mbox{\rm{Prob}} \left(  \zeta -  \beta \eta = ( b_i - \bfa_i^{\top} \bfv_1 ) \bfa^{\top} \bfv_2 \right)
      \, \mbox{d}\eta
      \end{eqnarray*}
      for all $\zeta \in \R$, which shows
      \begin{eqnarray*}
      \lefteqn{ \mbox{\rm{Prob}} \left( \bfa^{\top} \bfx = \beta \right) } \\
      & \leq &
       \int_{-\infty}^{\infty}
      \mbox{\rm{Prob}} \left( ( b_i - \bfa_i^{\top} \bfv_2 ) \bfa^{\top} \bfv_1 = \zeta \right) 
      + 
      \int_{-\infty}^{\infty}
      \mbox{\rm{Prob}} \left(  \zeta -  \beta \eta = ( b_i - \bfa_i^{\top} \bfv_1 ) \bfa^{\top} \bfv_2 \right)
      \, \mbox{d}\eta
      \, \mbox{d}\zeta .
      \end{eqnarray*}
      Now, for $\zeta \in \R$, $\zeta \neq 0$, 
      \begin{eqnarray*}
      \lefteqn{
      \mbox{\rm{Prob}} \left( ( b_i - \bfa_i^{\top} \bfv_2 ) \bfa^{\top} \bfv_1 = \zeta \right) 
      } \\
      & = & 
      \int_{-\infty}^{\infty}
      \mbox{\rm{Prob}} \left( ( b_i - \bfa_i^{\top} \bfv_2 ) \theta = \zeta \mbox{~and~} \bfa^{\top} \bfv_1 = \theta \right) 
      \, \mbox{d}\theta \\
      & \leq & 
      \int_{-\infty}^{\infty}
      \mbox{\rm{Prob}} \left( ( b_i - \bfa_i^{\top} \bfv_2 ) \theta = \zeta \right) 
      + \mbox{\rm{Prob}} \left( \bfa^{\top} \bfv_1 = \theta \right) 
      \, \mbox{d}\theta \\
      & = & 
      \int_{-\infty}^{\infty}
      \mbox{\rm{Prob}} \left( ( b_i - \bfa_i^{\top} \bfv_2 ) \theta = \zeta \right) 
      \, \mbox{d}\theta \\
      & = & 
      \int_{-\infty}^0
      \mbox{\rm{Prob}} \left( b_i - \bfa_i^{\top} \bfv_2 = \zeta/\theta \right) 
      \, \mbox{d}\theta
      +
      \int_{0}^{\infty}
      \mbox{\rm{Prob}} \left( b_i - \bfa_i^{\top} \bfv_2 = \zeta / \theta \right) 
      \, \mbox{d}\theta \\
      & = & 0 
      \end{eqnarray*}
      and 
      \begin{eqnarray*}
      \mbox{\rm{Prob}} \left( ( b_i - \bfa_i^{\top} \bfv_2 ) \bfa^{\top} \bfv_1 = 0 \right) 
      & = & 
      \mbox{\rm{Prob}} \left( b_i = \bfa_i^{\top} \bfv_2  \mbox{~or~} \bfa^{\top} \bfv_1 = 0 \right) \\
      & \leq & 
      \mbox{\rm{Prob}} \left( b_i = \bfa_i^{\top} \bfv_2 \right) + \mbox{\rm{Prob}} \left( \bfa^{\top} \bfv_1 = 0 \right) \\
      & = & 0 .
      \end{eqnarray*}
      In a similar fashion, it can be shown that
     \[
      \mbox{\rm{Prob}} \left(  \zeta -  \beta \eta = ( b_i - \bfa_i^{\top} \bfv_1 ) \bfa^{\top} \bfv_2 \right) = 0 
     \]
     for all $\zeta, \eta \in \R$.
      As a consequence, $\bfa^{\top} \bfx = \beta$ holds with probability
      zero. It is clear that the same analysis can be conducted for all other pairs
      of vectors $\bfv_{\ell_1}$, $\bfv_{\ell_2}$ with $\ell_1 < \ell_2$.
      $\Box$

\end{enumerate}

Lemma~\ref{lemma:ind}, invoked for $i = m$, shows part~1 of Theorem~\ref{theorem:main}.
It remains to discuss the complexity of the algorithm.  The 
\texttt{for}-loop is over $m$ steps, and each step involves three calls
to \textbf{rec}.  Executing \textbf{rec} costs two inner products of vectors
in~$\R^n$, two multiplications of vectors with scalars and one vector additon,
i.~e.\ the complexity of a call to \textbf{rec} is~$O(n)$. 
These considerations show Part~2 of Theorem~\ref{theorem:main}.

Some remarks are in order.

\begin{itemize}
\item It is clear that the algorithm also works for complex matrices 
      $\bfA \in \C^{m \times n}$ and complex right-hand sides $\bfb \in \C^n$.
      Again, no symmetry assumption on~$\bfA$ is necessary.
\item Some bookkeeping shows that the big-$O$ constant of the run time of
      the algorithm is ca.~$15$. While this appears large as compared to
      the big-$O$ constant of Gaussian elimination, 1/3, note that
      $15 n^2 < n^3 / 3$ for $n > 45$.
\item The algorithm is optimal in the sense that its run time is of the same
      order as its input size $(\bfA, \bfb)$.
\item The algorithm does not need to access the row vectors $\bfa_1, \ldots, \bfa_m$ directly;
      instead, it suffices to provide a routine that computes the action
      $\bfa_j^{\top} \bfv$ of a row $\bfa_j$ on an arbitrary vector~$\bfv \in \R^n$.
\item If the algorithm stops with failure in step~$k$, then we have
      $\bfa_j^{\top} \bfu_k = \bfa_j^{\top} \bfv_k = \bfa_j^{\top} \bfw_k = b_k$ for
      $j = 1, \ldots, k-1$, i.~e.\ the algorithm provides at least solutions
      to a subset of the system of equations.
\item Stability issues:  part of the stability of the algorithm rests on
      the size of quantities of the form $1 / ( \bfa_{k+1}^{\top} ( \bfu - \bfv ) )$.
      It is, at present, unclear how this quantity can be bounded away from zero.
\item In the exposition above, exactly $n+1$ vectors $\bfv_{\ell}$ are iteration
      vectors within the algorithm.  We can, of course, use more than $n+1$ 
      vectors to iterate over, and
      choose in each step $L > n+1$ pairs of vectors $\bfv_{i_1}, \bfv_{i_2}$ from
      the current iterates to feed into~\mbox{\textbf{rec}}. This increases
      the complexity of the algorithm from $O(n^2 m)$ to $O(L n m)$. However,
      choosing the right pairs of iterates $\bfv_{i_1}, \bfv_{i_2}$ in an adaptive
      fashion, possibly discarding results whose norm is too large, might
      alleviate the stability issues mentioned above.
\item Another way that might be useful to stabilize the method at hand is to 
      measure the degeneracy of a pair  $(\bfv_{i}, \bfv_{j})$ chosen in an
      iteration. If, say, $\Vert \bfv_i - \bfv_j \Vert$ is smaller than a certain
      threshold, the pair can either be discarded, or~$\bfv_j$ can be replaced by 
      $\bfv_i + c ( \bfv_j - \bfu_i )$ for
      a certain $c> 1$. A value of $0 < c < 1$ can be chosen if $\Vert \bfv_j - \bfv_i \Vert$
      grows too large.
\item In Step~3~(a), it is not necessary to always choose $n+1$ pairs of indices
      (and thus generate resp.\ update all of the vectors $\bfv_1, \ldots, \bfv_{n+1}$.
      Indeed, after $k$~steps of the main loop, all those vectors are in the $k$-dimensional
      affine subspace defined by the first $k$~equations $\bfa_i^{\top} \bfx = b_i$
      $(i=1, \ldots, k)$ and will remain in this subspace for all further iterations.
      Thus, after step~$k$, only $n+1-k$ pairs are needed to generate corresponding
      $n+1-k$ new vectors. 
\end{itemize}

\end{document}